\documentclass[12pt]{amsart}

\usepackage{amsfonts,amsmath,latexsym,amssymb,verbatim}
\usepackage{amsthm}

\theoremstyle{plain}
\newtheorem{THEOREM}{Theorem}[section]

\newtheorem{theorem}[THEOREM]{Theorem}

\newtheorem{lemma}[THEOREM]{Lemma}
\newtheorem{proposition}[THEOREM]{Proposition}

\theoremstyle{definition}

\newtheorem{definition}[THEOREM]{Definition}

\theoremstyle{remark}

\newtheorem{remark}[THEOREM]{Remark}

\newcommand{\N}{\ensuremath{\mathbb{N}}}   
\newcommand{\R}{\ensuremath{\mathbb{R}}}   
\newcommand{\T}{\ensuremath{\mathbb{T}}}   

\def \a {\alpha}
\def \b {\beta}
\def \d {\delta}

\def \e {\epsilon}
\def \f {\phi}

\def \l {\lambda}

\def \s {\sigma}

\def \t {\tau}

\def \be {\begin{equation}\label}

\def \ra {\rightarrow}
\def \ss {\subset}

\newcommand{\rest}[2]{#1\raisebox{-0.3ex}{\mbox{$\mid_{#2}$}}}

\def \be#1{\begin{equation}\label{#1}}
\def \besub#1{\begin{subequations}\label{#1}}
\def \endsub {\end{subequations}}

\def \emph#1{{\it #1}}
\def \textbf#1{{\bf #1}}

\begin{document}

\title{Nonlinear instability for the Navier-Stokes equations}

\author{Susan Friedlander \and Nata\v{s}a Pavlovi\'{c} \and Roman Shvydkoy}

\thanks{Friedlander received partial support from NSF grant DMS-0202767
and Pavlovic received partial support from NSF grant DMS-0304594.
Friedlander and Shvydkoy are very grateful to IAS for its
hospitality during the writing of this article.}
\address{University of Illinois, Department of Mathematics
Statistics and Computer Science, (m/c 249), Chicago, IL 60607 }
\email{susan@math.uic.edu; shvydkoy@math.uic.edu}
\address{Princeton
University, Department of Mathematics, Princeton,  NJ 08544}
\email{natasa@math.princeton.edu}


 \maketitle

\begin{abstract}
It is proved, using a bootstrap argument, that linear instability
implies nonlinear instability for the incompressible Navier-Stokes
equations in $L^p$ for all $p \in (1,\infty)$ and any finite or
infinite domain in any dimension $n$.
\end{abstract}

\section{Introduction}

The stability/instability of a flow of viscous incompressible fluid
governed by the Navier-Stokes equations is a classical subject with
a very extensive literature over more than $100$ years. Much of the
classical literature has concerned the stability of relatively
simple specific flows (e.g. Couette flows and Poiseuille flows), the
spectrum of the Navier-Stokes equations linearized about such flows
and the role of the critical Reynolds number delineating the
linearly stable and unstable regimes. An elegant result for general
bounded flows was proved by Serrin \cite{Serrin59} who used energy
methods to show that all flows are nonlinearly stable in $L^2$ norm
when the Reynolds number is less than a specific constant ($\pi
\sqrt{3}$). Hence all steady flows that are sufficiently slow or
sufficiently viscous are stable. However for many physical
situations the Reynolds number is much larger than $\pi \sqrt{3}$,
often by many orders of magnitude and observations indicate that
such flows are unstable.

Linear instability has been confirmed in some specific examples by
demonstrating existence of an nonempty unstable spectrum for the
linearized Navier-Stokes operator. For example, Meshalkin and Sinai
\cite{MS61} used Fourier series and continued fractions to show the
existence of unstable eigenvalues in the case of so called
Kolmogorov flows (i.e. plane parallel shear flow with a sinusoidal
profile). In a book published in Russian in 1984 (and in English in
1989) Yudovich \cite{Yudovich89} obtained an important result
relating linear stability/instability for the Navier-Stokes
equations with nonlinear stability/instability. These results were
proved in the function space $L^q(\Omega)$ with $q\geq n$ in
$n$-spatial dimensions. A fairly general abstract theorem of
Friedlander et al \cite{FSV97} can be applied to the Navier-Stokes
equations in a finite domain to prove nonlinear instability in
$H^s$, $s > \frac{n}{2}+1$ when the linearized operator has an
unstable eigenvalue in $L^2$.

In this present paper we extend the result that linear instability
implies nonlinear instability for the Navier-Stokes equations to all
$L^p$ spaces with $1 < p < \infty$ and both finite domains and
${\mathbb{R}}^{n}$. We note that our result includes nonlinear
instability in the $L^2$ energy norm which we claim is the natural
norm in which to consider issues of stability and instability.

The technique we employ to prove our main result is a bootstrap
argument. Such arguments have been previously employed by several
authors to prove under certain restrictions that linear instability
implies nonlinear instability for the 2 dimensional Euler equation
(Bardos et al \cite{BGS02}, Friedlander and Vishik \cite{FV03}, Lin
\cite{Lin}). Because in general the spectrum of the Euler operator
has a continuous component, unlike the Navier-Stokes operator in a
finite domain whose spectrum is purely discrete, these nonlinear
instability results for the Euler equation are much more limited
than those presented here for the Navier-Stokes equations.

\section{Notation and Formulation}

We consider solutions to the Navier-Stokes equations

\besub{NSE}
\begin{align}
 \frac{\partial q}{\partial t} &= - (q \cdot \nabla) q - \nabla p + {R}^{-1} \Delta q + f, \label{NS} \\
 \nabla \cdot q &= 0, \label{divfree}
\end{align}
\endsub
where $q(x,t)$ denotes the $n$-dimensional velocity vector, $p(x,t)$
denotes the pressure and $f(x)$ is an external force vector. The
dimensionless parameter $R$ is the Reynolds number defined as $R =
\frac{VL}{\nu}$ where $V$ and $L$ are characteristic velocity and
length scales of the system and $\nu$ is the viscosity of the fluid.
In Section \ref{finite} we consider the system on the
$n$-dimensional torus ${\mathbb T}^{n}$ and in a bounded domain
$\Omega \subset {\mathbb R}^{n}$. In Section \ref{infinite} we
consider the system in ${\mathbb R}^{n}$. The results are valid in
all dimensions $n$ although the most relevant physical cases are
$n=2$ and $3$. We impose the standard boundary conditions on
solutions of \eqref{NSE} for each type of domain. The results in
Sections 3 and 4 prove that spectral instability for the linearised
Navier-Stokes equations implies nonlinear instability in $L^p$ for
$1<p<\infty$. In Section 5 we prove a result relating spectral
stability with nonlinear stability in $L^p$ for $p>n$.

Here and thereafter, for any $p \in [1, \infty)$, $L^p$ denotes the
usual Lebesgue space, with norm denoted $\|\cdot\|_p$, intersected
with the space of divergence free functions. We let $W^{s,p}$ stand
for the Sobolev space in the same context with norm denoted
$\|\cdot\|_{s,p}$.

We consider an arbitrary steady solution of \eqref{NSE}
\begin{subequations}\label{eqns}
\begin{align}
0 &= - (U_0 \cdot \nabla) U_0 - \nabla P_0 + {R}^{-1} \Delta U_0 + f , \label{ns} \\
 \nabla \cdot U_0 &= 0. \label{eqdivfree}
\end{align}
\end{subequations}
We assume $U_0(x) \in C^{\infty}$ and $f(x) \in  C^{\infty}$. To
discuss stability of $U_0$ we rewrite the Navier-Stokes equations
\eqref{NSE} in perturbation form with $q(x,t) = U_0(x) + v(x,t)$
\besub{NSEpert}
\begin{align}
\frac{\partial v}{\partial t} & =  - (U_0 \cdot \nabla) v -  (v
\cdot \nabla) U_0 + {R}^{-1} \Delta v
-\nabla \cdot (v \otimes v) - \nabla p \label{perns} \\
\nabla \cdot v &= 0 \label{perdivfree} \\
\rest{v}{t=0} &= v_0
\end{align}
\endsub

Applying the Leray projector $\mathbb{P}$ onto the space of
divergence free functions, we write \eqref{perns} in the operator
form: \be{opns} \frac{\partial v}{\partial t} = A v + N(v,v)
\end{equation}
where
\begin{align}
Av &= {\mathbb{P}} [- (U_0 \cdot \nabla) v -  (v \cdot \nabla) U_0 + {R}^{-1} \Delta v] \label{A} \\
N(v,v) &=  {\mathbb{P}} [-\nabla \cdot (v \otimes v)] \label{nonlin}
\end{align}
We note that the linear operator $A$ is a bounded perturbation, to
lower order, of the Stokes operator $R^{-1}{\mathbb{P}} \Delta$. The
operator $A$ generates a strongly continuous semigroup in every
Sobolev space $W^{s,p}$ which we denote by $e^{At}$: \be{sem} v(t) =
e^{At} v_0, \; \; v_{0} \in W^{s,p}.
\end{equation}

We now define a suitable version of Lyapunov (nonlinear) stability
for the Navier-Stokes equations.

\begin{definition}\label{defstab} Let $(X,Z)$ be a pair of Banach
spaces. An equilibrium  $U_0$ which is the solution of \eqref{eqns}
is called $(X,Z)$ nonlinearly stable if, no matter how small $\rho >
0$, there exists $\delta > 0$ so that $v_0 \in X$ and \be{cond}
\|v_0\|_{Z} < \delta
\end{equation}
imply the following two assertions
\begin{enumerate}
\item[(i)] there exists a global in time solution to \eqref{NSEpert}
such that $v(t) \in C([0,\infty); X)$;
\item[(ii)] $\|v(t)\|_{Z} < \rho$ for a.e. $t \in [0,\infty)$.
\end{enumerate}
An equilibrium  $U_0$ that is not stable in the above sense is
called Lyapunov unstable.
\end{definition}

We will drop the reference to $(X,Z)$ where it does not lead to
confusion.

We note that under this strong definition of stability, loss of
existence of a solution to \eqref{opns} is a particular case of
instability. We remark that in literature there are many definitions
of a solution to the Navier-Stokes equations. These include
``classical'' solutions that are continuous functions of each
argument (and very few such solutions are known), ``weak'' solutions
defined via test functions by Leray \cite{Leray} and ``mild''
solutions introduced by  Kato-Fujita \cite{Kato-F62}. It is this
last concept of existence that we will invoke because we utilize a
``mild'' integral representation of the solution to \eqref{opns} via
Duhamel's formula. We remark that to date local in time existence of
mild solutions for the Navier-Stokes equations is proved only in
$L^{p}$, $p \geq n$, (for $p > n$ by Fabes-Jones-Riviere
\cite{FJR72} and for $p=n$ by Kato \cite{Kato84}). The existence of
weak solutions has been proved in $L^2$ by Leray \cite{Leray}, in
$L^p$ for all $2 \leq p < \infty$ by C. Calderon \cite{Calderon90},
and for uniformly locally square integrable initial data by
Lemari\'{e} \cite{Lemarie}. For a survey of existence results see
for example, Temam \cite{Temam2000} and Cannone \cite{Cannone04}.

We now state the main result of this paper:

\begin{theorem}\label{maintheorem}
Let $1<p<\infty$ be arbitrary. Suppose that the operator $A$ over
$L^p$ has spectrum in the right half of the complex plane. Then the
flow $U_0$ is $(L^q, L^p)$ nonlinearly unstable for any
$q>\max\{p,n\}$.
\end{theorem}
The proof of this theorem essentially uses properties of the
operator $A$ which are stated in Lemmas \ref{decaylemma} and
\ref{pertlemma}. The instability result is proved using a bootstrap
argument which is presented in Section \ref{finite} in the case of
finite domains $\Omega$ and ${\mathbb T}^{n}$ and in Section
\ref{infinite} in the case of ${\mathbb{R}}^{n}$.

Here we state a version of the Sobolev embedding theorem that we
shall invoke in the proof of Theorem \ref{maintheorem}.
\begin{proposition}\label{negSob}
Let $s>0$, $1<r_1<\infty$, and $1<r_2<\infty$ satisfy
\begin{equation}\label{Sobexp}
\frac{1}{r_1}  < 1 - \frac{s}{n}, \quad r_2  \leq r_1, \quad
\frac{1}{r_2} \leq \frac{1}{r_1} + \frac{s}{n}.
\end{equation}
Then \be{Sob-s-neg} \|f\|_{{-s, r_1}} \lesssim \|f\|_{r_2}
\end{equation}
\end{proposition}

\begin{proof}
Recall that for $s>0$ and $1<r<\infty$, $W^{-s,r}$ is defined as the
dual space to $W_0^{s,r'}$, where $1/r + 1/r' = 1$. The inequalities
\eqref{Sobexp} can be rewritten as
\begin{equation}\label{classSob}
sr_1'  < n, \quad r_1'  \leq r_2' \leq \frac{nr_1'}{n-sr_1'}.
\end{equation}
Thus, the standard Sobolev embedding theorem implies that
\be{stanSob} \|f\|_{{r_2'}} \lesssim \|f\|_{{s,r_1'}}.
\end{equation}
Applying \eqref{stanSob} we obtain
\begin{equation*}
\|f\|_{{-s,r_1}}  = \sup_{ \|g\|_{{s,r_1'}} \leq 1 } \langle
f,g\rangle
 \lesssim \sup_{ \|g\|_{{r_2'}} \leq 1 } \langle f,g \rangle  = \|f\|_{r_2},
\end{equation*}
which proves the proposition.
\end{proof}

\section{Finite domain} \label{finite}

In this section we present a proof of Theorem \ref{maintheorem} in
the case of finite domains ${\mathbb T}^{n}$ and $\Omega \ss \R^n$.

Let $\mu$ be the eigenvalue of $A$ with maximal positive real part,
which we denote by $\lambda$, and let $\phi \in L^p$, with
$\|\phi\|_p = 1$, be the corresponding eigenfunction. We note that
in the case of a finite domain all eigenfunctions of $A$ are
infinitely smooth.

For a fixed $0< \delta < \l$ we denote by $A_{\delta}$ the following
operator: \be{A_d} A_{\delta} = A- \lambda - \delta.
\end{equation}

Now we state two auxiliary lemmas which hold both in the case of a
finite and in the case of an infinite domain.
\begin{lemma}\label{decaylemma}
For every $0 < \alpha < 1$ and $p>1$ there exists a constant $M>0$
such that for all $t>0$ one has%
\be{analest} \|A_{\delta}^{\alpha} \; e^{A_{\delta}t} \|_{L^{p}
\rightarrow L^{p}} \leq \frac{M}{t^{\alpha}}.
\end{equation}
\end{lemma}
This lemma holds generally for any bounded analytic semigroup (see
\cite{Pazy}). The rescaling of $A$ given by \eqref{A_d} ensures that
the semigroup $e^{A_\d t}$ is bounded. The fact that it is analytic
is proved by Yudovich \cite{Yudovich89} and Giga \cite{Giga81}.

\begin{lemma}\label{pertlemma}
For every $1/2 < \alpha < 1$ and $p>1$ there exists a constant $C>0$
such that
\begin{equation}\label{giga}
\|A_{\delta}^{-\alpha} \; f\|_{{p}} \leq C \|f\|_{-2\alpha, p}.
\end{equation}
\end{lemma}
In the case of a bounded domain the lemma follows by duality from
the papers of Giga \cite{Giga-Stokes} and Seeley \cite{Seeley72}. On
the torus and $\R^n$ one can check \eqref{giga} directly using the
Fourier transform and integral representation for fractional power
of a generator \cite{Pazy}.

We are now in a position to prove Theorem \ref{maintheorem}.

Let us fix an arbitrary small $\e>0$, and solve the Cauchy problem
\eqref{NSEpert} with initial condition $v_0 = \e \f$. We note that
for such initial condition, with $\e$ small enough, there exists
unique global in time classical solution to \eqref{NSEpert} (see,
for example, \cite{Temam2000}). Using Duhamel's formula we write the
solution in the form

\begin{equation}\label{duhamel}
v(t) = \e e^{t\mu} \f + B(t),
\end{equation}
where
$$
B(t) =\int_{0}^{t} e^{A(t-\tau)} N(v,v)(\tau) \; d\tau.
$$

The main idea of the proof is to show that the bilinear term $B(t)$
grows at most like the square of the norm of $v(t)$ for as long as
the latter is bounded by a constant multiple of $\e e^{\l t}$. The
$L^q$-metric in which such control is possible has to satisfy the
assumption $q>n$. Since this condition is not assumed for $p$ we
will use $L^q$ as an auxiliary space, while our final instability
result will be proved  in $L^p$ as stated.

\begin{lemma}\label{L:best}
Let $q > n$. Then there exists a constant $C>0$ such that the
following estimate holds
\begin{equation}\label{E:best}
\|B(t)\|_q \leq C \int_0^t e^{(\l + \d)(t -
\t)}\frac{1}{(t-\t)^\a}\|v(\t)\|_q^2 \; d\t,
\end{equation}
for some $1/2 < \a < 1$.
\end{lemma}
\begin{proof}
Indeed, for any $0<\a<1$, we can write
$$
B(t) = \int_{0}^{t} e^{(\lambda +\delta)(t - \tau)} \;
A_{\delta}^{\alpha} \;  e^{A_{\delta}(t-\tau)}
{A_{\delta}^{-\alpha}} N(v,v)(\tau) \; d\tau.
$$
Hence, by Lemma \ref{decaylemma},
$$
\|B(t)\|_q \lesssim \int_{0}^{t} e^{(\lambda +\delta)(t - \tau)}
\frac{1}{(t-\t)^\a} \| {A_{\delta}^{-\alpha}} N(v,v)(\tau) \|_q\;
d\tau.
$$
By Lemma \ref{pertlemma}, we have
$$
\| {A_{\delta}^{-\alpha}} N(v,v) \|_q \lesssim \| N(v,v)\|_{-2\a,q}
\lesssim \|v\otimes v\|_{1-2\a,q},
$$
where the last inequality follows from the continuity of the Leray
projection. We now choose $\a$ sufficiently close to $1$ so that $
q> n/(2\a - 1)$. This would fulfill the conditions of Proposition
\ref{negSob} with $s = 2\a - 1$, $r_1 = q$ and $r_2 = q/2$. Thus,
\begin{equation}\label{replace}
\| {A_{\delta}^{-\alpha}} N(v,v) \|_q \lesssim \|v\otimes v\|_{q/2}
\lesssim \|v\|_q^2.
\end{equation}
Inserting this in the last estimate for $\|B(t)\|_q$ we finally
obtain \eqref{E:best}.
\end{proof}

\def \x {\mathfrak{X}}
Let us fix $q> \max\{n,p\}$. So, in particular, \eqref{E:best}
holds. For any $Q > \|\f\|_q$ let $T = T(Q)$ be the maximal time
such that
\begin{equation}\label{maxtime}
    \|v(t)\|_q \leq Q\e e^{\l t}, \quad \forall\; t\leq T.
\end{equation}
Notice that \eqref{maxtime} holds for $t = 0$. Hence, $T>0$ by
continuity. In fact, we show that this critical time $T$ is
sufficiently large for any choice of $Q$. First, let us observe that
for any $t \leq T$, by Lemma \ref{L:best},
$$
\|B(t)\|_q \leq C Q^2 \e^2 \int_0^t e^{(\l + \d)(t -
\t)}\frac{1}{(t-\t)^\a} e^{2\l \t} \; d\t.
$$
Splitting the integral into two integrals over $[0,t-1]$ and
$[t-1,t]$, one can show that it behaves asymptotically as $e^{2\l
t}$. Hence, perhaps with a different $C>0$ independent of $Q$ or
$t$, we obtain the following estimate
\begin{equation}\label{E:onb}
    \|B(t)\|_q \leq C(Q\e e^{\l t})^2, \quad \forall\; t \leq T.
\end{equation}
Using \eqref{E:onb} we now prove an estimate on the size of $T$.
\begin{lemma}\label{L:size}
For any $Q > \|\f\|_q$ one has the following inequality
\begin{equation}\label{E:size}
    \e e^{\l T} \geq \frac{Q - \|\f\|_q}{CQ^2}.
\end{equation}
\end{lemma}
\begin{proof}
If $T = \infty$, the inequality is trivial. If $T<\infty$, then at
time $t = T$ the inequality \eqref{maxtime} turns into equality and
we obtain using \eqref{duhamel} and \eqref{E:onb}
$$
Q \e e^{\l T} = \| v(T)\|_q \leq \e  e^{\l T} \|\f\|_q + C (Q\e
e^{\l T})^2.
$$
The lemma now easily follows.
\end{proof}
Let $\x_*$ denote the constant on the right hand side of
\eqref{E:size}. In view of \eqref{E:size} there exists  time $t_*
\leq T$ such that $\x_* = \e e^{\l t_*}$. Since $q > p$ we trivially
have
\begin{equation}\label{E:inclusion}
\|B(t)\|_p \leq C' \|B(t)\|_q,
\end{equation}
for some $C'>0$. So, by the triangle inequality applied to
\eqref{duhamel} we obtain
\begin{equation}\label{analogue}
\|v(t_*)\|_p \geq \x_* - C'C \x_*^2 = \x_*(1 - C'C \x_*).
\end{equation}
Since $C$ and $C'$ are independent of $Q$, we could choose $Q = Q_0$
in the beginning of the argument so close to $\|\f\|_q$  that $\x_*
< 1/(2C'C)$. Then
$$
\| v(t_*)\|_p \geq \x_*/2 = c_0.
$$

This finishes the proof of Theorem \ref{maintheorem} in the case of
a finite domain.

We remark that in the case of a finite domain our method proves a
stronger result. Since the eigenfunction $\f$ belongs to $C^\infty$,
the size of initial perturbation can be measured in the stronger
metric of $C^\infty$ so that $\|v_0\|_{C^\infty} \leq \e$, whereas
instability at the critical time $t_*$ is measured in the weak
$L^p$-metric.

\section{Infinite domain} \label{infinite}

The case of $\R^n$ brings two main difficulties to the proof. First,
we no  longer have the inclusion $L_q \subset L_p$ to satisfy
\eqref{E:inclusion}. Second, there may not be an exact smooth
eigenfunction $\phi$ corresponding to $\mu \in \sigma(A)$, because
the operator $A$ has a non-compact resolvent over $\R^n$.

\subsection{Estimates for $B(t)$}

In the case of $R^n$ we replace the single estimate
\eqref{E:inclusion} with a\ sequence of recursive estimates
improving integrability exponent on each step.

Let $L$ be the first integer such that $2^L p >n$. By Lemma
\ref{L:best}, which is valid on $R^n$ too, we have
\begin{equation}\label{nonlin1}
\|B(t)\|_{2^L p} \leq C \int_0^t e^{(\l + \d)(t -
\t)}\frac{1}{(t-\t)^\a}\|v(\t)\|_{2^Lp}^2 \; d\t,
\end{equation}
for some $1/2 < \a < 1$. On the other hand, for every
$l=0,\ldots,L-1$ one has, in place of \eqref{replace},
$$
\|A_d^{-\a}N(v,v)\|_{2^lp} \lesssim \|v \otimes v\|_{1-2\a,2^lp}
\lesssim \|v \otimes v\|_{2^lp} \lesssim \|v\|_{2^{l+1}p}^2.
$$
Thus, we obtain
\begin{equation}\label{nonlin2}
\|B(t)\|_{2^l p} \leq C \int_0^t e^{(\l + \d)(t -
\t)}\frac{1}{(t-\t)^\a}\|v(\t)\|_{2^{l+1}p}^2 \; d\t.
\end{equation}

We postpone the use of \eqref{nonlin1} and \eqref{nonlin2} till
Lemma \ref{L:key}, where we show the analogue of \eqref{analogue}
for the case of $R^n$.

\subsection{Construction of approximate eigenfunctions}

Suppose now that $\mu \in  \s(A)$ lies on the boundary of the
spectrum and has the greatest positive real part $\l$. In this case
there exists a sequence of functions $\{f_m\}_{m=1}^\infty \subset
L^p(\R^n)$ such that
\begin{gather*}
 \|f_m\|_{p} = 1,\\
 \lim_{m \ra \infty}\|Af_m - \mu f_m\|_{p} = 0,\\
\intertext{and as a consequence, for every $t>0$,} \lim_{m \ra
\infty}\|e^{tA}f_m - e^{t \mu} f_m\|_{p} = 0.
\end{gather*}

\begin{lemma}\label{L:appr}
There exists a sequence $\{\f_m\}_{m=1}^\infty \subset L^p(\R^n)$
such that the following is true
\begin{itemize}
\item[(i)] $\|\f_m\|_{p} = 1$, $n \in \N$;
\item[(ii)] For every $q > p $ there is a constant $M_q$ such that
$$
\|\f_m\|_{q} \leq M_q$$ for all $m \in \N$;
\item[(iii)] $\|e^{tA}\f_m\|_{p} \geq \frac{1}{2}e^{t\l}$, for all $0 \leq t \leq m$;
\item[(iv)] $\|e^{tA}\f_m\|_{q} \leq {2}\|\f_m\|_q e^{t\l}$, for all $0 \leq t \leq m$ and $p \leq q \leq 2^Lp$.
\end{itemize}
\end{lemma}

\begin{proof}
\def \tg {\tilde{\f}_m}
Let $\tg = e^Af_m$. Since
$$
\|e^Af_m - e^\mu f_m\|_{p} \ra 0
$$
we conclude that
\begin{equation}\label{gbound}
c \leq \|\tg\|_{p} \leq C,
\end{equation}
for all $m\in \N$. Denote $\f_m = \tg \cdot \|\tg\|_{p}^{-1}$.
Clearly, (i) is satisfied. To prove the other three statements we
fix $s>0$ such that $n = sp$. By the end-point Sobolev embedding
theorem and \eqref{gbound} we have, for any $q>p$,
\begin{align*}
\|\f_m\|_q & \lesssim \|\tg\|_{q} = \|e^Af_m\|_{q} \lesssim
\|e^Af_m\|_{{s,p}} \\&\lesssim \|A_\delta^s e^A f_m\|_{p} \lesssim
\|f_m\|_{p} = 1.
\end{align*}
This proves (ii).

Furthermore, we have
\begin{align*}
\|e^{tA}\f_m - e^{t \mu }\f_m\|_{q} & \lesssim \|e^A(e^{tA}f_m - e^{t \mu }f_m)\|_{{s,p}}
\lesssim \|A_\d^s e^A(e^{tA}f_m - e^{t \mu }f_m)\|_{p}\\
&\lesssim \|e^{tA}f_m - e^{t \mu}f_m\|_{p} \ra 0,
\end{align*}
as $m \ra \infty$ for each fixed $t>0$ and $p\leq q$. So, by
choosing an appropriate subsequence, we achieve (iii) and (iv).
\end{proof}

\subsection{Bootstrap argument}
Let us fix an arbitrary $\e >0$ and find $m \in \N$ such that
\begin{equation}\label{m}
    \e e^{\l m} >1.
\end{equation}
This $m$ will be fixed though the rest of the argument. We solve the
Cauchy problem \eqref{NSEpert} with initial condition $v_0 = \e
\f_m$. Lemma \ref{L:appr} shows that $\f_m \in L^q$ uniformly in $m$
for all $q>p$. In particular, for any fixed $q > \max\{p,n\}$ there
exists a mild solution in  $Z = L^q$ for which the Duhamel
formulation holds:
\begin{equation}\label{duh}
  v(t) = \e e^{At} \f_m + B(t).
\end{equation}
We note that failure for $v(t)$ to satisfy \eqref{duh} for all $t>0$
or being in $C([0,\infty),X)$ is regarded as instability by
definition. We thus can assume in the rest of the argument that
\eqref{duh} holds for all $t>0$ and $v \in C([0,\infty),X)$. In
addition, since $\f_m \in L^{2^Lp}$ and $2^Lp>n$, the solution
$v(t)$ belongs to $L^{2^Lp}$ at least for a certain initial period
of time. Our subsequent estimates will show that, in fact, $v(t) \in
L^{2^Lp}$ over a time interval of the order $\log {1/\e}$.

Let $Q > 2 \|\f_m\|_{2^L p}$ be arbitrary, and define $T = T(Q)$ to
be the maximal time such that
\begin{equation}\label{max}
    \|v(t)\|_{2^L p} \leq Q \e e^{\l t}, \text{ for all } t \leq T.
\end{equation}
Like in the previous section the following inequality holds
\begin{equation}\label{infinest}
    \|B(t)\|_{2^L p} \leq C (Q\e e^{\l t})^2, \quad \forall t \leq
    T.
\end{equation}

\begin{lemma}\label{L:estT}
For any $Q > 2 \|\f_m\|_{2^L p}$ we have
\begin{equation}\label{E:estT}
\e e^{\l T} \geq \min\left\{1; \frac{ Q - 2\|\f_m\|_{2^L p} }{ C
Q^2} \right\},
\end{equation}
where $C>0$ is independent of $Q$.
\end{lemma}
\begin{proof}
If $T \geq m$, we appeal to \eqref{m}. If $T < m$, then at $t = T$
the inequality \eqref{max} must turn into equality. Thus, in view of
\eqref{infinest} we have
\begin{equation*}
Q \e e^{\l T} = \|v(T(Q))\|_{2^L p } \leq 2 \|\f_m\|_{2^L p} \e
e^{\l T} + C (Q\e e^{\l T})^2,
\end{equation*}
which implies \eqref{E:estT}.
\end{proof}

We will choose $Q$ appropriately after the following key lemma.

\begin{lemma}\label{L:key}
There are constants $C_2,\ldots,C_{2^{L+1}}$ and $2\leq K\leq
2^{L+1}$ independent of $Q$ and $m$ such that for any $t \leq
\min\{T, m\}$ one has the following inequality
\begin{multline}\label{E:key}
\|v(t)\|_p \geq \frac{1}{2}\x - C_2 \x^2 - \ldots - C_{K-1} \x^{K-1} -  \\
-C_K Q^K \x^K - \ldots  - C_{2^{L+1}} Q^{ 2^{L+1}} \x^{2^{L+1}},
\end{multline}
where $\x = \e e^{\l t}$.
\end{lemma}

\begin{proof}
First we bound all the norms $\|v(t)\|_{2^l p}$, $l=1,\ldots,L$ from
above using the estimates on the nonlinear term \eqref{nonlin2},
\eqref{infinest}. We start with $l=L$ and invoke \eqref{infinest} to
obtain
$$
\|v(t)\|_{2^L p} \leq \e e^{\l t} \|\f_m\|_{2^L p} + C (\e e^{\l t})
^2 Q^2 = C_1 \x + C_2 Q^2 \x^2,
$$
for all $t \leq T$. We note that our constants may change during the
proof.

By the previous inequality and \eqref{nonlin2} with $l = L-1$ we obtain%
\begin{align*}
\|v(t)\|_{2^{L-1} p} & \leq \x  \|\f_m\|_{2^{L-1} p} + \\
&+ \int_0^t  \frac{e^{(t - \t)(\l + \d)}}{(t - \t)^{\a}}(C_1 \e
e^{\l \t} + C_2
Q^2 \e^2 e^{2 \l \t} )^2 d \t \\
& \leq C_1 \x + C_2 \x^2 + C_3 Q^3 \x^3 + C_4 Q^4 \x^4.
\end{align*}
Here and thereafter we use the fact that for any $k \geq 2$ one has
$$
\int_0^t e^{(t-\t)(\l + \d)} (t-\t)^{-\a} e^{k\l \t} d\t \lesssim
e^{k \l t}.
$$

By induction on $l$ we arrive at
$$
\|v(t) \|_{2p} \leq C_1 \x + \ldots + C_{\tilde{K}-1} \x^{\tilde{K}
- 1} +C_{\tilde{K}} Q^{\tilde{K}}\x^{\tilde{K}} + \ldots + C_{2^L}
Q^{2^L} \x^{2^L},
$$
and hence,
$$
\|B(t)\|_p \leq C_2 \x^2 + \ldots + C_{K-1} \x^{K-1} + C_K Q^K \x^K
+ \ldots  + C_{2^{L+1}} Q^{ 2^{L+1}} \x^{2^{L+1}}.
$$
Finally, using (iii) of Lemma \eqref{L:appr} and the triangle
inequality on \eqref{duh} in the opposite direction to get
\eqref{E:key}.
\end{proof}

We will choose a $Q$ so that the RHS of \eqref{E:key} is bigger than
an absolute constant at
$$
\x=\x_* = \min\left\{1; \frac{ Q - 2\|\f_m\|_{2^L p} }{ C Q^2}
\right\}.
$$
For this $\x_*$, due to Lemma \ref{L:estT} and our initial
assumption \eqref{m}, there exists a $t_* \leq \min\{T,m\}$ such
that $\x_* = \e e^{t_* \l}$. Hence, Lemma \ref{L:key} applies to
obtain instability at time $t=t_*$.

It is convenient to seek $Q$ in the form
$$
Q = (2 + a \|\f_m\|_{2^L p} ) \|\f_m\|_{2^L p},
$$
where $0<a<1$. Then
\begin{equation*}
\frac{ Q - 2\|\f_m\|_{2^L p} }{ C Q^2} = \frac{a}{C(2 + a
\|\f_m\|_{2^L p})^2} \leq \frac{a}{4C}.
\end{equation*}
Choosing $a < 4C$ we ensure that
$$
\frac{ Q - 2\|\f_m\|_{2^L p} }{ C Q^2} <1
$$
and hence,
$$
\x_* =\frac{ Q - 2\|\f_m\|_{2^L p} }{ C Q^2}.
$$
By the above estimate and (ii) of Lemma \ref{L:appr}, we have
$$
\frac{a}{C(2+M_{2^Lp})^2} \leq \x_* \leq \frac{a}{4C},
$$
or
\begin{equation}\label{exest}
    \frac{a}{C'} \leq \x_* \leq \frac{a}{c'}.
\end{equation}
We notice that since $Q$ is bounded by a constant independent of $m$
and $a$, we can bound the minimum
$$
\min\left\{ 1 ; \min_{2\leq k \leq K-1} (C_k 4^k)^{-1/k} ; \min_{K
\leq k \leq 2^{L+1}} (C_k Q^k 4^k)^{-1/k} \right\}
$$
from below by some constant $c_0$ independent of $Q$. Let $a =
\min\{4C,c'c_0/2\}$. Then from \eqref{exest}, we obtain
$$
\tilde{c}_0 \leq \x_* \leq c_0.
$$
Thus, by \eqref{E:key},
$$
\|v(t_*)\|_p \geq \tilde{c}_0(\frac{1}{2} - \frac{1}{16} - \ldots )
= c.
$$

This finishes the proof.

We remark again that like in the case of a finite domain our method
yields a slightly stronger result. Since $\f_m \in W^{s,p}$
uniformly, we can measure the size of initial perturbation in the
metric of any Sobolev space $W^{s,p}$ for all $s>0$.

\section{Stability result}
Bootstrap techniques can also be used to prove that linear stability
implies nonlinear stability for the Navier-Stokes equations in $L^q$
for $q>n$. In particular this reproves the classical stability
theorem of Yudovich \cite{Yudovich89}.

\begin{theorem}
Let $q>n$ be arbitrary. Assume the operator $A$ in $L^q$ has
spectrum confined to the left half of the complex plane. Then the
flow $U_0$ is $(L^q,L^q)$ nonlinearly stable. The result holds in
$\T^n$ and $\Omega$, and in any spatial dimension $n$.
\end{theorem}
\begin{proof}
We recall that that any analytic semigroup possesses the spectral
mapping property. From the assumption that the spectrum of $A$ is
confined to the left half plane we thus conclude that the
exponential type of the semigroup $e^{At}$ is negative. Hence, there
exists $\l >0$ such that
\begin{equation}\label{expdecay}
    \|e^{At} v_0\|_q \leq M e^{-\l t} \|v_0\|_q,
\end{equation}
for all $t>0$ and $v_0 \in L^q$. From Duhamel's formula
\eqref{duhamel} with the initial condition replaced by $v_0$, and by
argument similar to that used in the proof of Lemma \ref{L:best}, we
have
\begin{equation}\label{stabest}
\|v(t)\|_q \leq M e^{-\l t} \|v_0\|_q + C \int_0^t e^{-\l(t-\t)}
(t-\t)^{-\a} \|v(\t)\|_q^2 d \t.
\end{equation}

Again let $T$ be the maximal time for which
\begin{equation}\label{maxstab}
    \|v(t)\|_q \leq 2M \|v_0\|_q e^{-\l t}, \quad t \leq T.
\end{equation}
Combining \eqref{stabest} and \eqref{maxstab} gives
\begin{align*}
\|v(t)\|_q & \leq M e^{-\l t} \|v_0\|_q + 4 M^2 C e^{-2\l t}
\|v_0\|_q^2  \\
& \leq M e^{-\l t} \|v_0\|_q (1 + 4M C\|v_0\|_q),
\end{align*}
for $t\leq T$. We choose $\|v_0\|_q < (8MC)^{-1}$. Then the previous
inequality implies that
\begin{equation}\label{smaller}
\|v(t)\|_q \leq \frac{3}{2}M\|v_0\|_q e^{-\l t},
\end{equation}
for $t \leq T$. Hence, the assumption of \eqref{maxstab} implies the
smaller bound of \eqref{smaller}, which gives a contradiction with a
maximal finite $T$. Thus, $T = \infty$ and the bound \eqref{maxstab}
holds for all $t \geq 0$. This bound implies the global existence of
the solution to \eqref{NSE} and condition (ii) of Definition
\ref{defstab} for a sufficiently small choice of $\|v_0\|_q$.

\end{proof}

\begin{remark} The instability/stability results in this paper can
be generalized to all the equations of motion that are augmented
versions of the equations for incompressible, dissipative fluids
described in operator form by an appropriate version of
\eqref{opns}. This includes the magnetohydrodynamic equations for a
dissipative electrically conducting fluid, the equations for an
incompressible, stratified fluid with viscous and thermal
dissipation and the so called modified Navier-Stokes equations with
$(-\Delta)$ replaced by $(-\Delta)^\b$ where $\b > 1/2$.
\end{remark}


\def\cprime{$'$} \def\cprime{$'$} \def\cprime{$'$} \def\cprime{$'$}
  \def\cprime{$'$} \def\cprime{$'$} \def\cprime{$'$} \def\cprime{$'$}
  \def\cprime{$'$}
\providecommand{\bysame}{\leavevmode\hbox
to3em{\hrulefill}\thinspace}

\end{document}